\documentclass[a4paper,11pt]{amsart}

\usepackage{a4wide}
\usepackage{amssymb}
\usepackage{mathrsfs}
\usepackage{mathtools}
\usepackage{enumerate}
\usepackage{titletoc}
\usepackage[colorlinks=true, urlcolor=blue, linkcolor=blue, citecolor=green]{hyperref}
\usepackage{mathscinet}

\theoremstyle{plain}
\newtheorem{theorem}{Theorem}[section]
\newtheorem{lemma}[theorem]{Lemma}
\newtheorem{corollary}[theorem]{Corollary}
\newtheorem{proposition}[theorem]{Proposition}
\newtheorem*{conjecture*}{Conjecture}
\newtheorem*{question*}{Question}

\theoremstyle{definition}
\newtheorem{definition}[theorem]{Definition}
\newtheorem{remark}[theorem]{Remark}

\newtheorem*{notation}{Notation}

\numberwithin{equation}{section}

\makeatletter
\@namedef{subjclassname@2020}{\textup{2020} Mathematics Subject Classification}
\makeatother

\title[Time derivative estimates]{Time derivative estimates for normalized parabolic $p$-Laplace equations}

\author{Se-Chan Lee}
\address{School of Mathematics, Korea Institute for Advanced Study, Seoul 02455, Republic of Korea}
\email{sechan@kias.re.kr}

\author{Hyungsung Yun}
\address{Department of Mathematics, Changwon National University, 20 Changwondaehak-ro, Uichang-gu, Changwon-si, Gyeongsangnam-do, 51140, Republic of Korea}
\email{hyungsung@changwon.ac.kr}  



\thanks{Se-Chan Lee has been supported by the KIAS Individual Grant (No. MG099001) at Korea Institute for Advanced Study. Hyungsung Yun has been supported by the National Research Foundation of Korea (NRF) grant funded by the Korea government (MSIT) (RS-2026-25469633).}

\begin{document}

\begin{abstract}
We establish the local boundedness of time derivatives of bounded viscosity solutions to the normalized parabolic $p$-Laplace equation for every $p\in(1,\infty)$. The proof combines a logarithmic Bernstein argument with quantitative stability to derive a recursive estimate between different regularizations. For $p=1$ and $n\geq2$, we construct a bounded viscosity solution which is not H\"older continuous in time with any positive exponent.
\end{abstract}

\maketitle

\section{Introduction}
In this paper, we are concerned with the interior Lipschitz regularity in the time variable of bounded viscosity solutions $u$ to the normalized parabolic $p$-Laplace equation
\begin{equation}\label{eq:main}
 u_t=\Delta_p^N u
 \coloneqq |Du|^{2-p}\Delta_pu
 =\Delta u+(p-2)\frac{u_i u_j}{|Du|^2}u_{ij}
\end{equation}
in $Q_1$, where $p\in(1,\infty)$ and $Q_r\coloneqq B_r\times(-r^2,0]$ for $r>0$. 

Tug-of-war games provide game-theoretic interpretations of elliptic equations involving the infinity-Laplacian in \cite{PSSW09} and, with noise, the normalized $p$-Laplacian for $p\in(1,\infty)$ in \cite{PS08}. Asymptotic mean value characterizations of $p$-harmonic
functions and viscosity solutions of normalized parabolic $p$-Laplace equations were established in \cite{MPR10a} and \cite{MPR10b}, respectively. Time-dependent games and their connections with normalized parabolic equations were further studied in \cite{Han22, PR16}. We refer to \cite{Doe11} for applications of the corresponding evolution equations to image processing. At $p=1$, equation \eqref{eq:main} reduces to the level-set mean curvature flow equation; see \cite{CGG91,ES91}.

If $p=2$, the equation \eqref{eq:main} is the heat equation, and interior estimates for solutions are well known. For general $p\in(1,\infty)$, its coefficient matrix is
\begin{equation*}
    a_0^{ij}(q)\coloneqq\delta_{ij}+(p-2)\frac{q_iq_j}{|q|^2} \quad\text{for $q \in \mathbb{R}^n \setminus\{0\}$}.
\end{equation*}
Since the eigenvalues of $(a_0^{ij})$ lie between $\min\{1,p-1\}$ and $\max\{1,p-1\}$, the equation is uniformly parabolic for every $p \in (1, \infty)$. However, when $n\geq2$ and $p\ne2$, the coefficient matrix is discontinuous where the gradient vanishes. This singular dependence on the gradient prevents a direct application of the classical regularity theory for equations with smooth coefficients.

The regularity theory for \eqref{eq:main} has been developed in several directions. Jin--Silvestre~\cite{JS17} established interior H\"older gradient estimates for \eqref{eq:main}, while the corresponding inhomogeneous equation $u_t=\Delta_p^Nu+f$ was treated in \cite{ABM23, AP18}. In particular, the regularity result in \cite{JS17} implies that
\begin{equation*}
    u \in C_{t, \mathrm{loc}}^{\alpha}(Q_1) \quad \text{ for some $\alpha=\alpha(n, p)\in(0,1)$.}
\end{equation*}
Extensions of these results to singular and degenerate parabolic equations 
\begin{equation}\label{eq-general-main}
    u_t=|Du|^{\gamma}\Delta_p^Nu \quad \text{for $p \in (1, \infty)$ and $\gamma \in (-1, \infty)$}
\end{equation}
were obtained in \cite{IJS19}; see also \cite{Att20,AR20} for the inhomogeneous case. Moreover, Fang--Zhang established gradient H\"older estimates for normalized equations with variable exponent $p(x,t)$ in \cite{FZ21} and for equations with degeneracy or singularity of double-phase type in \cite{FZ23}. We also refer to \cite{AS22,LY24} for perturbation-type regularity results for normalized equations, to \cite{LLY24} for degenerate and singular fully nonlinear equations, and to \cite{APR17,BM20} for the corresponding elliptic theory. 

The local integrability of weak time derivatives of solutions to \eqref{eq:main} has also been studied. H{\o}eg--Lindqvist~\cite{HL20} proved that $u_t\in L^2_{\mathrm{loc}}(Q_1)$ for $1<p<14/5$. Dong--Peng--Zhang--Zhou~\cite{DPZZ20} improved both the range of $p$ and the integrability exponent by establishing 
\begin{equation*}
u_t\in L^q_{\mathrm{loc}}(Q_1) \quad \text{for some $q=q(n,p)>2$,}    
\end{equation*}
for all $1<p<\infty$ when $n=2$ and for $1<p<3+2/(n-2)$ when $n\geq3$. For \eqref{eq-general-main} with suitable ranges of $p$ and $\gamma$, Feng--Parviainen--Sarsa~\cite{FPS23,FPS25} obtained local $L^2$ estimates for $u_t$ when $n=2$. 

For sufficiently smooth data, Lipschitz regularity in time for solutions to \eqref{eq:main} was established for Cauchy--Neumann problems~\cite[Theorem~4.1]{Doe11} and Cauchy--Dirichlet problems on suitable domains~\cite{BG13, BG15}. The constants in these estimates depend on $C^2$ bounds for the prescribed data.

\medskip

In this paper, we prove an interior $L^\infty$ estimate for $u_t$ when $p\in(1,\infty)$, with a bound depending only on $n$, $p$ and $\|u\|_{L^\infty(Q_1)}$. This provides a positive answer to the question raised by H{\o}eg--Lindqvist~\cite{HL20} of whether $u_t$ is locally bounded without smoothness assumptions on the lateral boundary data.

\begin{theorem}\label{thm:main}
Let $u\in C(Q_1)$ be a bounded viscosity solution to \eqref{eq:main} in $Q_1$ with $p\in(1,\infty)$. Then $u$ is locally Lipschitz continuous in time. In particular, $u$ is weakly differentiable in time and $u_t\in L^\infty_{\mathrm{loc}}(Q_1)$ with the uniform estimate
\begin{equation*}
 \|u_t\|_{L^\infty(Q_{1/2})}\leq C\|u\|_{L^\infty(Q_1)},
\end{equation*}
where $C>0$ is a constant depending only on $n$ and $p$.
\end{theorem}

Related interior $L^\infty$ estimates for $u_t$ were obtained by Lee--Lian--Yun--Zhang~\cite{LLYZ25} for \eqref{eq-general-main} with $p\in(1,\infty)$ and $\gamma \in (0, \infty)$, while the present work deals with the critical case $\gamma=0$. As in that work, the time derivative estimate in Theorem~\ref{thm:main} allows us to deduce improved spatial regularity by reduction to an elliptic problem. Since this estimate does not provide continuity of $u_t$, we use an elliptic resolvent approximation to obtain continuous right-hand sides with uniform bounds, and then apply \cite[Theorem~1.3]{APR17}.

Let $\alpha_p=\alpha_p(n)\in(0,1]$ denote the supremum of the exponents $\alpha\in(0,1]$ for which uniform interior $C^{1,\alpha}$ estimates hold for $p$-harmonic functions.
We obtain the following consequence of Theorem~\ref{thm:main}.

\begin{corollary}[Nearly optimal gradient regularity]\label{cor:almost-optimal-gradient}
Let $u\in C(Q_1)$ be a bounded viscosity solution to \eqref{eq:main} in $Q_1$ with $p\in(1,\infty)$. For every $\beta\in(0,\alpha_p)$, we have $u \in C^{1,\beta}_{\mathrm{loc}}({Q}_{1})$ with the uniform estimate
\begin{equation*}
\|u\|_{C^{1,\beta}(\overline{Q}_{1/2})} \leq C\|u\|_{L^\infty(Q_1)},
\end{equation*}
where $C>0$ depends only on $n$, $p$ and $\beta$. 
\end{corollary}

Corollary~\ref{cor:almost-optimal-gradient} refines the gradient H\"older estimate of Jin--Silvestre~\cite{JS17}, which holds for some exponent $\alpha=\alpha(n,p)>0$, by allowing every exponent $\beta<\alpha_p$. The spatial range is nearly optimal, since stationary solutions of \eqref{eq:main} include $p$-harmonic functions.

The restriction $p>1$ in Theorem~\ref{thm:main} cannot be removed. Setting $p=1$ in \eqref{eq:main} gives the level-set mean curvature flow equation
\begin{equation}\label{eq:mean-curvature}
 u_t=\Delta u-\frac{u_i u_j}{|Du|^2}u_{ij}.
\end{equation}
We construct the following example by applying a continuous nondecreasing function to a smooth solution, exploiting the invariance property of the equation~\cite{CGG91,Gig06}. 

\begin{theorem}\label{thm:p-one}
Let $n\geq2$. There exists a bounded viscosity solution to \eqref{eq:mean-curvature} in $Q_1$ which is not H\"older continuous in time with any positive exponent at an interior point.
\end{theorem}

\begin{remark}
The other limiting case is the normalized infinity heat equation
\begin{equation*}
u_t=\Delta_\infty^N u \coloneqq \frac{\langle D^2uDu,Du\rangle}{|Du|^2}.
\end{equation*}
This equation is the formal limit of $u_t=p^{-1}\Delta_p^N u$ as $p\to\infty$. Juutinen--Kawohl~\cite{JK06} proved interior Lipschitz estimates in space, which, together with a standard barrier argument, imply local $1/2$-H\"older continuity of $u$ in time. Related estimates for equations with a bias parameter and a source term were studied by Liu--Jiang~\cite{LJ19}. Our interior time Lipschitz estimate concerns $p \in (1, \infty)$ and does not cover this limiting equation.
\end{remark}

Let us illustrate the main ideas in the proof of Theorem~\ref{thm:main}. Following the Bernstein argument in \cite{LLYZ25}, we first derive a preliminary time derivative estimate (Lemma~\ref{lem:logarithmic}) for the regularized solutions $u^\varepsilon$ constructed in Proposition~\ref{prop:regularized}. In the critical case $\gamma=0$, however, the power-type auxiliary function used in \cite{LLYZ25} is not appropriate for deriving this estimate. Instead, we introduce a logarithmic-type auxiliary function, which allows us to control the bad terms in the linearized equation but gives an estimate depending on $\varepsilon$ in a logarithmic way.

To obtain a uniform estimate in $\varepsilon$, we compare different regularizations by distinguishing two regimes according to the size of the gradient. In the small gradient regime, the sharper bound for the auxiliary function gives an additional small factor in the error term. The maximum principle then allows us to reduce the problem to estimating the time derivative where the gradient is bounded away from zero. In the large gradient regime, we combine the interior Schauder estimates with the quantitative stability estimate of Kurkinen--Liu~\cite{KL26} to compare the time derivative with that of a regularized solution with a larger parameter.

Combining the estimates in the two regimes, we obtain a recursive estimate between different regularizations on nested cylinders. More precisely, writing $l=\log(e/\varepsilon)$, for $l$
sufficiently large we can choose a larger regularization parameter $\delta>\varepsilon$ such that
\begin{equation*}
    \|u_t^\varepsilon\|_{L^\infty(Q_s)}
    \leq \|u_t^\delta\|_{L^\infty(Q_{s+2/l})}
    +C/l \quad \text{and} \quad
    \log(e/\delta)\leq l/2;
\end{equation*}
see Section~\ref{sec:time-estimates} for details. We iterate this comparison until the logarithmic parameter falls below a fixed threshold. Since the accumulated errors are uniformly bounded and all cylinders remain inside a fixed cylinder, the preliminary Bernstein estimate at the final step guarantees a bound independent of $\varepsilon$. By passing to the limit, we finish the proof of Theorem~\ref{thm:main}.
\medskip

The paper is organized as follows. In Section~\ref{sec:preliminaries}, we introduce viscosity solutions, construct the regularized solutions, and prove quantitative stability estimates. Section~\ref{sec:log-lipschitz} establishes the logarithmic time derivative estimate, based on the Bernstein technique. In Section~\ref{sec:time-estimates}, we derive the estimates in two different regimes and combine them with an iteration argument to prove Theorem~\ref{thm:main}. Finally, in Section~\ref{sec:p-one}, we construct the example in Theorem~\ref{thm:p-one}.

\section{Preliminaries}\label{sec:preliminaries}

In this section, we summarize the notation, definitions, and regularity results that will be used throughout the paper.

\begin{notation}
We denote by $\mathcal{S}^n$ the space of real symmetric $n\times n$ matrices and use the Frobenius norm
\begin{equation*}
 \|M\|\coloneqq\left(\sum_{i,j=1}^n|M_{ij}|^2\right)^{1/2} \quad \text{for $M=(M_{ij})$}.
\end{equation*}
For $r>0$ and $(x_0,t_0)\in\mathbb R^n\times\mathbb R$, we set
\begin{equation*}
 Q_r(x_0,t_0)\coloneqq B_r(x_0)\times(t_0-r^2,t_0] \quad \text{and} \quad Q_r=Q_r(0,0). 
\end{equation*}
The parabolic boundary of $Q_r(x_0,t_0)$ is
\begin{equation*}
 \partial_pQ_r(x_0,t_0)
 \coloneqq\bigl(\overline{B}_r(x_0)\times\{t_0-r^2\}\bigr)
 \cup\bigl(\partial B_r(x_0)\times[t_0-r^2,t_0]\bigr).
\end{equation*}
We write $Du$ and $D^2u$ for the spatial gradient and Hessian, respectively, and sum repeated spatial indices from $1$ to $n$. 
For a cylinder $Q$ and $\alpha\in(0,1)$, we use the following definitions for parabolic H\"older spaces; see \cite{IS13book}:
\begin{itemize}
 \item We say $u\in C^\alpha(\overline{Q})$ if $u$ is $\alpha$-H\"older continuous in $x$ and $\alpha/2$-H\"older continuous in $t$. We define the $C^{\alpha}$-norm of $u$ by
\begin{equation*}
    \|u\|_{C^\alpha(\overline{Q})}\coloneqq \|u\|_{L^\infty(Q)} + [u]_{C_x^\alpha(\overline{Q})} + [u]_{C_t^{\alpha/2}(\overline{Q})}.
\end{equation*}
 \item We say $u\in C^{1,\alpha}(\overline{Q})$ if $Du \in C^{\alpha}(\overline Q)$ and $u$ is $(1+\alpha)/2$-H\"older continuous in $t$. We define the $C^{1,\alpha}$-norm of $u$ by
\begin{equation*}
    \|u\|_{C^{1,\alpha}(\overline{Q})}\coloneqq \|u\|_{L^\infty(Q)} + \|Du\|_{C^{\alpha}(\overline{Q})} + [u]_{C_t^{(1+\alpha)/2}(\overline{Q})}.
\end{equation*}
\end{itemize}
\end{notation}

We use the standard viscosity definition based on the upper and lower semicontinuous envelopes of the operator; see \cite{CIL92}. For the normalized parabolic $p$-Laplace equation, this definition appears in \cite[Definition~2.1]{MPR10b} and \cite[Definition~2.1]{AP18}. 

To be precise, set for $q\in\mathbb{R}^n\setminus\{0\}$ and $X\in\mathcal{S}^n$,
\begin{equation*}
 F(q,X)\coloneqq\operatorname{tr} X+(p-2)\frac{\langle Xq,q\rangle}{|q|^2}.
\end{equation*}
Its upper and lower semicontinuous envelopes are defined by
\begin{equation*}
 F^*(q,X)\coloneqq
 \limsup_{\substack{(z,Y)\to(q,X)\\z\ne0}}F(z,Y)
 \quad \text{and} \quad
 F_*(q,X)\coloneqq
 \liminf_{\substack{(z,Y)\to(q,X)\\z\ne0}}F(z,Y).
\end{equation*}
For $q\ne0$, continuity gives 
\begin{equation*}
    F^*(q,X)=F_*(q,X)=F(q,X).     
\end{equation*}
At $q=0$, they are given by
\begin{equation*}
 \begin{aligned}
 F^*(0,X)&=
 \begin{cases}
  \operatorname{tr} X+(p-2)\lambda_{\max}(X) &  \text{if $p\geq2$},\\
  \operatorname{tr} X+(p-2)\lambda_{\min}(X)& \text{if $1<p<2$},
 \end{cases}\\
 F_*(0,X)&=
 \begin{cases}
  \operatorname{tr} X+(p-2)\lambda_{\min}(X)&\text{if $p\geq2$},\\
  \operatorname{tr} X+(p-2)\lambda_{\max}(X)&\text{if $1<p<2$},
 \end{cases}
 \end{aligned}
\end{equation*}
where $\lambda_{\min}(X)$ and $\lambda_{\max}(X)$ are the smallest and largest eigenvalues of $X$, respectively.

\begin{definition}[Viscosity solutions]
Let $\mathcal O\subset\mathbb R^n\times\mathbb R$ be an open set. An upper semicontinuous function $u:\mathcal O\to\mathbb R$ is a viscosity subsolution to \eqref{eq:main} in $\mathcal O$ if, for every $\varphi\in C^2(\mathcal O)$ and every point
$(x_0,t_0)\in\mathcal O$ at which $u-\varphi$ attains a local maximum, we have
\begin{equation*}
\varphi_t(x_0,t_0) \leq F^*(D\varphi(x_0,t_0),D^2\varphi(x_0,t_0)).
\end{equation*}
A lower semicontinuous function $u:\mathcal O\to\mathbb R$ is a viscosity supersolution to \eqref{eq:main} in $\mathcal O$ if, for every $\varphi\in C^2(\mathcal O)$ and every point
$(x_0,t_0)\in\mathcal O$ at which $u-\varphi$ attains a local minimum, we have
\begin{equation*}
\varphi_t(x_0,t_0) \geq F_*(D\varphi(x_0,t_0),D^2\varphi(x_0,t_0)).
\end{equation*}
A function $u\in C(\mathcal O)$ is a viscosity solution to \eqref{eq:main} in $\mathcal O$ if it is both a viscosity subsolution and a viscosity supersolution.
\end{definition}

We introduce the regularized coefficient matrix
\begin{equation}\label{eq:matrix}
 a^{ij}_\varepsilon(q)\coloneqq \delta_{ij}+(p-2)\frac{q_iq_j}{\varepsilon^2+|q|^2}
 \quad \text{for $\varepsilon \in (0, 1]$ and $q \in \mathbb{R}^n$}.
\end{equation}
If we set $\lambda_p\coloneqq\min\{1,p-1\}$ and $\Lambda_p\coloneqq\max\{1,p-1\}$, then we have
\begin{equation}\label{eq:ellipticity}
 \lambda_p|\xi|^2
 \leq a_\varepsilon^{ij}(q)\xi_i\xi_j
 \leq\Lambda_p|\xi|^2
 \quad \text{for every $q,\xi\in\mathbb{R}^n$}.
\end{equation}
In particular, the regularized matrix $(a^{ij}_{\varepsilon})$ is uniformly elliptic. 

\begin{proposition}\label{prop:regularized}
Let $u \in C(Q_1)$ be a viscosity solution to \eqref{eq:main} in $Q_1$ with $\|u\|_{L^\infty(Q_1)}\leq1$. For every $\varepsilon\in(0,1]$, the Cauchy--Dirichlet problem
\begin{equation*}
 \left\{\begin{aligned}
  u_t^\varepsilon&=a_\varepsilon^{ij}(Du^\varepsilon)u_{ij}^\varepsilon &&\text{in }Q_{7/8},\\
  u^\varepsilon&=u &&\text{on }\partial_pQ_{7/8}
 \end{aligned}\right.
\end{equation*}
has a unique solution $u^\varepsilon\in C^\infty(Q_{7/8})\cap C(\overline{Q}_{7/8})$. Moreover, there exist $\alpha\in(0,1)$ and $C>0$, depending only on $n$ and $p$, such that
\begin{equation}\label{eq:gradient}
    \|u^\varepsilon\|_{C^{1,\alpha}(\overline{Q}_{3/4})}\leq C \quad \text{and}  \quad \|u^\varepsilon\|_{C^\alpha(\overline{Q}_{7/8})}\leq C.
\end{equation}
\end{proposition}

\begin{proof}
For each fixed $\varepsilon>0$, existence and uniqueness follow from the classical theory of quasilinear parabolic equations; see \cite[Chapter~VI, Theorem~4.4]{LSU68} and \cite[Lemma~5.1]{JS17} for details. The first interior estimate in \eqref{eq:gradient} is obtained from \cite[Theorem~4.5]{JS17}, after scaling and a standard covering argument.

To prove the second global estimate in \eqref{eq:gradient}, we let $\mathcal P^\pm$ denote the Pucci extremal operators with ellipticity constants $\lambda_p$, $\Lambda_p$, namely
\begin{equation*}
 \mathcal P^-(X)=\inf_{\lambda_p I\leq A\leq\Lambda_p I}\operatorname{tr}(AX)
 \quad \text{and} \quad
 \mathcal P^+(X)=\sup_{\lambda_p I\leq A\leq\Lambda_p I}\operatorname{tr}(AX).
\end{equation*}
By \eqref{eq:ellipticity} and the definition of viscosity solutions, both $w=u$ and $w=u^\varepsilon$ satisfy
\begin{equation*}
 \mathcal P^-(D^2w)\leq w_t\leq\mathcal P^+(D^2w)
 \quad\text{in $Q_{7/8}$}
\end{equation*}
in the viscosity sense. Since the interior H\"older estimate for $u$ first gives $\|u\|_{C^{\alpha_1}(\overline{Q}_{15/16})}\leq C$ for some $\alpha_1=\alpha_1(n,p)\in(0,1)$, the common boundary data $u|_{\partial_pQ_{7/8}}$ on $\partial_pQ_{7/8}$ have a uniformly bounded H\"older norm. Therefore, the boundary H\"older estimate for uniformly parabolic equations, combined with the interior estimate, now yields
\begin{equation*}
 \|u^\varepsilon\|_{C^{\alpha_2}(\overline{Q}_{7/8})}\leq C
\end{equation*}
for some $\alpha_2=\alpha_2(n,p)\in(0,\alpha_1]$; see \cite{Lie96, Wan92a,Wan92b} for instance. The estimate depends only on $n$ and $p$; in particular, it is uniform in $\varepsilon$. Choosing $\alpha$ smaller and $C$ larger, if necessary, finishes the proof of \eqref{eq:gradient}.
\end{proof}

The following quantitative stability estimate is essentially contained in \cite[Theorem~4.7(1)]{KL26}. For completeness, we give a direct proof based on their comparison argument, using the viscosity definition above and covering the full range $\varepsilon\in[0,1]$.

\begin{proposition}[Quantitative stability]\label{prop:quantitative-stability}
Let $p\in(1,\infty)$, let $\Omega\subset\mathbb R^n$ be a bounded domain, and set $\Omega_T=\Omega\times(0,T)$ for $T>0$. For $\varepsilon\in[0,1]$, let
$u^\varepsilon\in C(\overline{\Omega}_T)$ be a viscosity solution to
\begin{equation*}
\left\{
\begin{aligned}
u_t^\varepsilon&=a_\varepsilon^{ij}(Du^\varepsilon)u_{ij}^\varepsilon &&\text{in }\Omega_T,\\
u^\varepsilon&=u^0 &&\text{on }\partial_p\Omega_T.
\end{aligned}
\right.
\end{equation*}
 Suppose that there exist $C_0>0$ and $\theta\in(0,1]$ such that
\begin{equation*}
    |u^\varepsilon(x_1,t)-u^\varepsilon(x_2,t)| \leq C_0|x_1-x_2|^\theta \quad
    \text{for all $x_1,x_2\in\overline\Omega$, $t\in[0,T]$ and $\varepsilon\in[0,1]$.}
\end{equation*}
Then
\begin{equation*}
\|u^\varepsilon-u^0\|_{L^\infty(\Omega_T)} \leq C\varepsilon^{\theta/4} \quad\text{for all }\varepsilon\in(0,1],
\end{equation*}
where $C>0$ depends only on $n$, $p$, $T$, $\theta$ and $C_0$.
\end{proposition}

\begin{proof}
We use the matrix notation $a_\varepsilon=(a_\varepsilon^{ij})$
and fix $\varepsilon\in(0,1]$. For $q \neq 0$, it follows from a direct calculation that
\begin{equation*}
\begin{aligned}
\|a_\varepsilon(q)^{1/2}-a_0(q)^{1/2}\|^2
&=\left(\sqrt{1+(p-2)\frac{|q|^2}{\varepsilon^2+|q|^2}}-\sqrt{p-1}\right)^2\\
&\leq |p-2|\frac{\varepsilon^2}{\varepsilon^2+|q|^2}.
\end{aligned}
\end{equation*}

For an arbitrary $b>0$ and $M>0$ chosen later, define $\Phi: \overline\Omega\times\overline\Omega\times[0,T) \to \mathbb{R}$ by
\begin{equation*}
\Phi(x,y,t)=u^\varepsilon(x,t)-u^0(y,t)-\frac{|x-y|^4}{4\varepsilon}-M\varepsilon^{\theta/4}t-\frac{b}{T-t}.
\end{equation*}
Let $\Phi$ attain its maximum at some $(\bar x,\bar y,\bar t)$, and set $r=|\bar x-\bar y|$. Comparing this maximum with $\Phi(\bar x,\bar x,\bar t)$
and using the spatial H\"older estimate of $u^0$, we obtain
\begin{equation*}
\frac{r^4}{4\varepsilon}
\leq u^0(\bar x,\bar t)-u^0(\bar y,\bar t)\leq C_0r^\theta,
\end{equation*}
and so there exists $C>0$ depending only on $\theta$ and $C_0$ such that
\begin{equation*}
r\leq C\varepsilon^{1/(4-\theta)}
\leq C\varepsilon^{1/4}.
\end{equation*}
If $\bar t=0$ or if either $\bar x$ or $\bar y$ belongs to $\partial\Omega$, then the common boundary values and the spatial H\"older estimate give
\begin{equation*}
\max \Phi=\Phi(\bar x, \bar y, \bar t) \leq u^\varepsilon(\bar x,\bar t)-u^0(\bar y,\bar t) \leq C_0r^\theta \leq C\varepsilon^{\theta/4}.
\end{equation*}

We next claim that, for $M$ sufficiently large depending only on $n$ and $p$, the maximum cannot occur with $\bar x,\bar y\in\Omega$ and $\bar t\in(0,T)$. Suppose otherwise. By the parabolic version of Ishii--Jensen's lemma~\cite[Theorem~8.3]{CIL92},
there exist $\tau_1,\tau_2\in\mathbb R$
and $X,Y\in\mathcal S^n$ such that
\begin{flalign*}
&\text{(i)}\quad
\tau_1\leq a_\varepsilon^{ij}(q)X_{ij} \quad \text{and} \quad  \tau_2\geq F_*(q,Y) \quad \text{for }q=\frac{r^2(\bar x-\bar y)}{\varepsilon};&&\\
&\text{(ii)}\quad
\begin{pmatrix}
X&0\\
0&-Y
\end{pmatrix}
\leq
\frac{6r^2}{\varepsilon}
\begin{pmatrix}
I&-I\\
-I&I
\end{pmatrix};&&\\
&\text{(iii)}\quad
\tau_1-\tau_2=M\varepsilon^{\theta/4}+\frac{b}{(T-\bar t)^2}.&&
\end{flalign*}
If $r=0$, then $q=0$ and $X\leq 0 \leq Y$ from (ii). Therefore, the viscosity inequalities in (i) yield
\begin{equation*}
\tau_1\leq\operatorname{tr}X\leq0 \quad \text{and} \quad \tau_2\geq F_*(0,Y)\geq 0,
\end{equation*}
which contradicts $\tau_1-\tau_2>0$.

Suppose now that $r>0$. Applying the matrix inequality (ii) to $(a_\varepsilon(q)^{1/2}e_i,a_0(q)^{1/2}e_i)$,
where $e_i$ are the standard basis vectors in $\mathbb{R}^n$,
we observe that
\begin{equation*}
    e_i^Ta_{\varepsilon}(q)^{1/2}Xa_{\varepsilon}(q)^{1/2}e_i-e_i^Ta_0(q)^{1/2}Ya_{0}(q)^{1/2}e_i \leq \frac{6r^2}{\varepsilon}\left|(a_{\varepsilon}(q)^{1/2}-a_{0}(q)^{1/2})e_i\right|^2,
\end{equation*}
and so
\begin{equation*}
    \operatorname{tr}(a_{\varepsilon}(q)X)-\operatorname{tr}(a_0(q)Y) \leq \frac{6r^2}{\varepsilon}\|a_{\varepsilon}(q)^{1/2}-a_{0}(q)^{1/2}\|^2.
\end{equation*}
We then use (i) and (iii) to arrive at
\begin{equation*}
\begin{aligned}
M\varepsilon^{\theta/4}\leq \tau_1-\tau_2&\leq
\operatorname{tr}(a_\varepsilon(q)X)-\operatorname{tr}(a_0(q)Y)\\
&\leq
\frac{6r^2}{\varepsilon}\|a_\varepsilon(q)^{1/2}-a_0(q)^{1/2}\|^2\\
&\leq C\frac{r^2}{\varepsilon}\frac{\varepsilon^2}{\varepsilon^2+|q|^2}=C\frac{\varepsilon^3r^2}{\varepsilon^4+r^6},
\end{aligned}
\end{equation*}
where we used $|q|=r^3/\varepsilon>0$. Moreover, by considering two cases $r\leq\varepsilon^{2/3}$ and $r>\varepsilon^{2/3}$, we conclude that
\begin{equation*}
    M\varepsilon^{\theta/4} \leq C\frac{\varepsilon^3r^2}{\varepsilon^4+r^6}\leq C\varepsilon^{1/3}\leq C\varepsilon^{\theta/4}.
\end{equation*}
Taking $M>C$ gives a contradiction and proves the claim.

In conclusion, the maximum of $\Phi$ is attained when $\bar t=0$ or when either $\bar x$ or $\bar y$ belongs to $\partial \Omega$; in particular, $\max\Phi\leq C\varepsilon^{\theta/4}$. Evaluating $\Phi$ at $(x,x,t)$ and then letting $b\to0$,
we find
\begin{equation*}
u^\varepsilon(x,t)-u^0(x,t) \leq (C+MT)\varepsilon^{\theta/4} \quad\text{for }(x,t)\in\Omega_T.
\end{equation*}
Applying the same argument to $-u^\varepsilon$ and $-u^0$ gives the reverse estimate and completes the proof.
\end{proof}

We finally combine Propositions~\ref{prop:regularized} and \ref{prop:quantitative-stability} to compare the regularized solutions with $u$ and with each other.

\begin{corollary}\label{cor:approximation}
Let $\alpha \in (0,1)$ be chosen in Proposition~\ref{prop:regularized}. Then there exists $C=C(n,p)>0$ such that
\begin{equation}\label{eq:stability}
 \|u^\varepsilon-u\|_{L^\infty(Q_{7/8})}\leq C\varepsilon^{\alpha/4}
 \quad \text{for $\varepsilon \in (0,1]$}.
\end{equation}
In particular, if $0<\varepsilon \leq \delta \leq 1$, then we have
\begin{equation}\label{eq:pairwise}
 \|u^\varepsilon-u^\delta\|_{L^\infty(Q_{7/8})}\leq 2C\delta^{\alpha/4}.
\end{equation}
\end{corollary}

\begin{proof}
By \eqref{eq:gradient} and the interior H\"older estimate for $u$ established in the proof of Proposition~\ref{prop:regularized}, the solutions $u^\varepsilon$ and $u$ satisfy the assumptions of Proposition~\ref{prop:quantitative-stability} on $Q_{7/8}$, after translating the time variable. Thus, \eqref{eq:stability} and \eqref{eq:pairwise} hold with the exponent $\alpha/4$.
\end{proof}

\section{A logarithmic Bernstein estimate}\label{sec:log-lipschitz}
Throughout this and the next section, we fix the exponent $\alpha$ from Proposition~\ref{prop:regularized}, and assume that 
\begin{itemize}
    \item $u \in C(Q_1)$ is a bounded viscosity solution to \eqref{eq:main} in $Q_1$ with $\|u\|_{L^\infty(Q_1)}\leq1$;

    \item for $\varepsilon \in (0,1]$, $u^{\varepsilon} \in C^{\infty}(Q_{7/8}) \cap C(\overline{Q}_{7/8})$ is the regularized solution that satisfies the Cauchy--Dirichlet problem given in Proposition~\ref{prop:regularized}.
\end{itemize}

We begin with a preliminary time derivative estimate for the regularized solutions, with a bound depending on $\varepsilon$. The proof follows the Bernstein argument in \cite{LLYZ25}, with a logarithmic auxiliary function in place of a power of $\varepsilon^2+|Du^\varepsilon|^2$.

\begin{lemma}[Logarithmic time derivative estimates]\label{lem:logarithmic}
Let $1/2\leq s<s+d\leq2/3$ and $\varepsilon \in (0,1]$. Then we have
\begin{equation}\label{eq:logarithmic}
 \|u_t^\varepsilon\|_{L^\infty(Q_s)}
 \leq Cd^{-2}\log\frac e\varepsilon,
\end{equation}
where $C>0$ depends only on $n$ and $p$.
\end{lemma}

\begin{proof}
For notational convenience, we omit $\varepsilon$ from $u^\varepsilon$ and $a_\varepsilon^{ij}$. We consider the linearized operator
\begin{equation*}
 Lw\coloneqq a^{ij}(Du)w_{ij}+b^l w_l \quad \text{for $b^l\coloneqq a^{ij}_{q_l}(Du)u_{ij}$}.
\end{equation*}
We also write $b=(b^1, \ldots, b^n)$. Differentiation of the equation $u_t=a^{ij}(Du)u_{ij}$ with respect to the variables $t$ and $x_k$ gives, respectively:
\begin{equation}\label{eq:linearized}
 (\partial_t-L)u_t=0
 \quad \text{and} \quad
 (\partial_t-L)u_k=0.
\end{equation}
Moreover, it follows from \eqref{eq:matrix} and \eqref{eq:ellipticity} that
\begin{equation}\label{eq:drift}
    \begin{aligned}
    |u_t|&=|a^{ij}(Du)u_{ij}|\leq C\|D^2u\|,\\
    |b^l|&=|p-2|\left|\left[\frac{\delta_{il}u_j+\delta_{jl}u_i}{\varepsilon^2+|Du|^2}-\frac{2u_iu_ju_l}{(\varepsilon^2+|Du|^2)^2}\right] u_{ij} \right|\leq\frac{C\|D^2u\|}{(\varepsilon^2+|Du|^2)^{1/2}}.
    \end{aligned}
\end{equation}

For the smooth convex function
\begin{equation*}
\Psi_\varepsilon(q)\coloneqq\rho\log\frac\rho\varepsilon \quad \text{for $\rho=(\varepsilon^2+|q|^2)^{1/2}$},
\end{equation*}
it follows from a direct calculation that
\begin{equation*}
    D^2_q\Psi_{\varepsilon}(q)=\frac{1+\log(\rho/\varepsilon)}{\rho}I-\frac{\log(\rho/\varepsilon)}{\rho^3}q \otimes q \geq \frac{1}{\rho}I
\end{equation*}
and that
\begin{equation*}
    \begin{aligned}
        (\partial_t-L)\Psi_\varepsilon(Du)&=(\Psi_{\varepsilon})_{q_k}(Du) (\partial_t u_k-a^{ij}(Du)u_{kij}-b^lu_{kl}) \\
        &\quad-a^{ij}(Du)(\Psi_\varepsilon)_{q_kq_l}(Du)u_{ki}u_{lj}\\
        &=-a^{ij}(Du)(\Psi_\varepsilon)_{q_kq_l}(Du)u_{ki}u_{lj}.
    \end{aligned}
\end{equation*}
Since $(a^{ij}) \geq \lambda_pI$ and $D^2\Psi_{\varepsilon} \geq \rho^{-1}I$, we have
\begin{equation*}   
        K \coloneqq a^{ij}(Du)(\Psi_\varepsilon)_{q_kq_l}(Du)u_{ki}u_{lj}\geq\frac{\lambda_p\|D^2u\|^2}{(\varepsilon^2+|Du|^2)^{1/2}},
\end{equation*}
and so \eqref{eq:drift} yields that
\begin{equation}\label{eq:absorb}
 |u_t||b|\leq CK.
\end{equation}
On the other hand, the gradient bound \eqref{eq:gradient} and $\varepsilon \in (0,1]$ imply
\begin{equation*}
     0\leq\Psi_\varepsilon(Du)\leq C\log\frac e\varepsilon \quad \text{and} \quad 
     K\geq c\|D^2u\|^2 \quad\text{in $Q_{3/4}$}.
\end{equation*}

We now choose a smooth cut-off function $\eta$ on $\overline{Q}_{s+d}$ such that $0\leq\eta\leq1$, $\eta=1$ on $Q_s$, and $\eta=0$ on $\partial_pQ_{s+d}$, with
\begin{equation*}
    |\eta_t|+\|D^2\eta\|\leq Cd^{-2}, \quad
    |D\eta|\leq Cd^{-1}\quad\text{and}\quad
    |D\eta|^2\leq Cd^{-2}\eta.
\end{equation*}
We next let
\begin{equation*}
 A\coloneqq\|\eta u_t\|_{L^\infty(Q_{s+d})}.
\end{equation*}
There is nothing to prove if $A=0$, so we assume that $A>0$. It follows from \eqref{eq:linearized} that
\begin{equation*}
    \begin{aligned}
        (\partial_t-L)(\eta^2u_t^2)=&2\eta(\eta_t-a^{ij}\eta_{ij}-b^l\eta_l)u_t^2
 -2a^{ij}\eta_i\eta_j u_t^2\\
 &-2\eta^2a^{ij}u_{ti}u_{tj}-8\eta u_ta^{ij}\eta_i u_{tj}\\
 &\leq 2\eta(|\eta_t|+|a^{ij}\eta_{ij}|+|b||D\eta|)u_t^2
 -\eta^2a^{ij}u_{ti}u_{tj}+8|\eta u_ta^{ij}\eta_i u_{tj}|.
    \end{aligned}
\end{equation*}
By Young's inequality together with $(a^{ij}) \leq \Lambda_pI$, we have
\begin{equation*}
 8|\eta u_ta^{ij}\eta_i u_{tj}| \leq 8|\eta u_t|(a^{ij}\eta_i\eta_j)^{1/2}(a^{ij}u_{ti}u_{tj})^{1/2}\leq \eta^2a^{ij}u_{ti}u_{tj}+C|D\eta|^2u_t^2.
\end{equation*}
Moreover, we observe that
\begin{itemize}
    \item since $|\eta_t|+|a^{ij}\eta_{ij}| \leq Cd^{-2}$, we have
    \begin{equation*}
        2\eta(|\eta_t|+|a^{ij}\eta_{ij}|)u_t^2 \leq Cd^{-2} \eta u_t^2 \leq Cd^{-2}A|u_t|;
    \end{equation*}

    \item since $|D\eta| \leq Cd^{-1}$, we have
    \begin{equation*}
        2\eta |b||D\eta|u_t^2 \leq Cd^{-1} \eta u_t^2 |b| \leq Cd^{-1} A |u_t| |b|;
    \end{equation*}

    \item since $|D\eta|^2 \leq Cd^{-2}\eta$, we have
    \begin{equation*}
        C|D\eta|^2u_t^2 \leq Cd^{-2}\eta u_t^2 \leq Cd^{-2}A|u_t|.
    \end{equation*}
\end{itemize}
A combination of all estimates gives that
\begin{equation}\label{eq:product}
 (\partial_t-L)(\eta^2u_t^2)
 \leq Cd^{-2}A\|D^2u\|+Cd^{-1}A|u_t||b|.
\end{equation}
Thus, \eqref{eq:absorb} and $\|D^2u\|\leq (\|D^2u\|^2+1)/2\leq CK+1/2$ give
\begin{equation*}
    (\partial_t-L)(\eta^2u_t^2) \leq Cd^{-2}A(K+1)=-Cd^{-2}A (\partial_t-L)(\Psi_{\varepsilon}(Du)-t).
\end{equation*}
We then choose a sufficiently large constant $C_1>0$, depending only on $n$ and $p$, so that the auxiliary function
\begin{equation*}
     v\coloneqq\eta^2u_t^2 + C_1d^{-2}A(\Psi_\varepsilon(Du)-t)
\end{equation*}
satisfies $(\partial_t-L)v<0$ in $Q_{s+d}$. Since
\begin{equation*}
     v=C_1d^{-2}A\Psi_\varepsilon(Du)-C_1d^{-2}At\leq Cd^{-2}A\log\frac e\varepsilon \quad \text{on $\partial_pQ_{s+d}$},
\end{equation*}
the maximum principle gives the same bound for $v$ in $Q_{s+d}$. Finally, we arrive at
\begin{equation*}
    A^2=\|\eta u_t\|_{L^{\infty}(Q_{s+d})}^2 \leq Cd^{-2}A\log\frac e\varepsilon,
\end{equation*}
which finishes the proof by recalling that $\|u_t\|_{L^{\infty}(Q_s)} \leq A$.
\end{proof}

\begin{remark}[Log-Lipschitz regularity in time]\label{rem:log-lipschitz}
Lemma~\ref{lem:logarithmic} and Corollary~\ref{cor:approximation} already imply that $u$ is locally log-Lipschitz continuous in time. Indeed, taking $s=1/2$ and $d=1/16$ in \eqref{eq:logarithmic}, we obtain
\begin{equation*}
 |u(x,t+h)-u(x,t)|
 \leq 2\|u-u^\varepsilon\|_{L^\infty(Q_{7/8})}
      +h\|u_t^\varepsilon\|_{L^\infty(Q_{1/2})}
 \leq C\varepsilon^{\alpha/4}+Ch\log\frac e\varepsilon
\end{equation*}
whenever $x\in B_{1/2}$, $h>0$ and $t,t+h\in(-1/4,0]$. Choosing $\varepsilon=h^{4/\alpha}$ gives
\begin{equation*}
 |u(x,t+h)-u(x,t)|\leq Ch\log\frac {e}{h}.
\end{equation*}
\end{remark}

\section{Lipschitz regularity in time}\label{sec:time-estimates}
The goal of this section is to improve the log-Lipschitz estimate in Remark~\ref{rem:log-lipschitz} to the Lipschitz estimate in Theorem~\ref{thm:main}. To this end, we compare different regularizations by distinguishing the two regimes described in the introduction. 
\begin{itemize}
    \item In the small gradient regime $\{|Du^\varepsilon|<r\}$, the Bernstein auxiliary function has a smaller bound, and the maximum principle reduces the estimate to the set $\{|Du^\varepsilon|\geq r\}$, as established in Lemma~\ref{lem:small-gradient}. 
    \item In the large gradient regime $\{|Du^\varepsilon| \geq r\}$, local Schauder estimates combined with the quantitative stability estimate, allow us to compare time difference quotients at different regularization parameters, yielding the estimate in Lemma~\ref{lem:noncritical}.
\end{itemize}

For $1/2\leq s\leq2/3$ and $r>0$, we set
\begin{equation*}
\begin{aligned}
    M_\varepsilon(s)&\coloneqq\sup\left\{|u_t^\varepsilon(x,t)|:(x,t)\in Q_s\right\} \,(=\|u_t^\varepsilon\|_{L^\infty(Q_s)}),\\
    \widetilde M_\varepsilon(s,r)&\coloneqq \sup\left\{|u_t^\varepsilon(x,t)|: (x,t)\in Q_s,\ |Du^\varepsilon(x,t)|\geq r\right\}.
 \end{aligned}
\end{equation*}
Here we set $\widetilde M_\varepsilon(s,r)=0$ when the set $\{(x,t)\in Q_s:|Du^\varepsilon(x,t)|\geq r\}$ is empty. Both quantities are finite for each fixed $\varepsilon>0$, and $0\leq\widetilde M_\varepsilon(s,r)\leq M_\varepsilon(s)$. We first compare these two quantities and then estimate $\widetilde M_\varepsilon(s,r)$ by investigating a time derivative bound for a regularized solution with a larger parameter. 

\begin{lemma}\label{lem:small-gradient}
Let $1/2\leq s<s+d\leq2/3$ and $0< \varepsilon\leq r$. Then we have
\begin{equation}\label{eq:small-gradient}
    M_\varepsilon(s)\leq\widetilde M_\varepsilon(s+d,r)+Cd^{-2}r\log\left(e+\frac r\varepsilon\right),
\end{equation}
where $C>0$ depends only on $n$ and $p$.
\end{lemma}

\begin{proof}
We essentially repeat the proof of Lemma~\ref{lem:logarithmic}, again writing $u=u^\varepsilon$ for simplicity. Indeed, let $L$, $\Psi_\varepsilon$, $K$, $\eta$ be as in that proof, and set $A=\|\eta u_t\|_{L^\infty(Q_{s+d})}>0$.

On $\{|Du|\leq r\} \cap Q_{s+d}$, we can replace the bound $\rho=\sqrt{\varepsilon^2+|Du|^2} \leq C$ by $\rho \leq \sqrt{2}r$. This sharper bound together with the assumption $\varepsilon\leq r$ gives 
\begin{equation*}
     0\leq\Psi_\varepsilon(Du) \leq Cr\log\left(e+\frac{r}{\varepsilon}\right), \quad 
     K \geq cr^{-1}\|D^2u\|^2\quad\text{and}\quad 
     \|D^2u\| \leq CK+Cr.
\end{equation*}
Thus, \eqref{eq:absorb} and \eqref{eq:product} imply
\begin{equation*}
    (\partial_t-L)(\eta^2u_t^2)\leq Cd^{-2}A(K+r)\quad\text{on }\{|Du|<r\}\cap Q_{s+d}.
\end{equation*}
Therefore, we choose a sufficiently large $C_1$, depending only on $n$ and $p$, so that
\begin{equation*}
    v\coloneqq\eta^2u_t^2+C_1d^{-2}A(\Psi_\varepsilon(Du)-rt)
\end{equation*}
satisfies $(\partial_t-L)v<0$ on this set.

On $\{|Du|=r\} \cap Q_{s+d}$, the first term of $v$ is at most $A\widetilde M_\varepsilon(s+d,r)$, while it vanishes on $\partial_pQ_{s+d}$. Since the remaining terms are bounded by $Cd^{-2}Ar\log(e+r/\varepsilon)$, the maximum principle shows that
\begin{equation*}
    v\leq A\widetilde M_\varepsilon(s+d,r)+Cd^{-2}Ar\log\left(e+\frac r\varepsilon\right)
 \quad\text{in $\{|Du|<r\}\cap Q_{s+d}$}.
\end{equation*}
Indeed, if this bound failed, the maximum of $v$ on $\overline{Q}_{s+d}\cap\{|Du|\leq r\}$ would occur where $|Du|<r$, away from the initial and lateral boundary. At such a point, we have 
\begin{equation*}
    Dv=0, \quad D^2v\leq 0 \quad\text{and}\quad v_t\geq0,     
\end{equation*}
which contradicts $(\partial_t-L)v<0$. 

We are now ready to prove the estimate \eqref{eq:small-gradient}. If $\eta^2u_t^2$ attains its maximum at $(x_1, t_1) \in \{|Du|\geq r\} \cap Q_{s+d}$, then the definition of $\widetilde{M}_{\varepsilon}$ implies 
\begin{equation*}
    A=\eta(x_1, t_1)|u_t(x_1, t_1)| \leq |u_t(x_1, t_1)| \leq \widetilde{M}_{\varepsilon}(s+d, r).
\end{equation*}
Otherwise, $v\geq\eta^2u_t^2$ implies
\begin{equation*}
    A^2\leq A\widetilde M_\varepsilon(s+d,r)+Cd^{-2}Ar\log\left(e+\frac r\varepsilon\right).
\end{equation*}
Therefore, we have in either case
\begin{equation*}
 A\leq\widetilde M_\varepsilon(s+d,r)+Cd^{-2}r\log\left(e+\frac r\varepsilon\right).\qedhere
\end{equation*}
\end{proof}

We next estimate the time derivative where the gradient is bounded away from zero. We first use interior Schauder estimates to approximate $u_t^\varepsilon$ by a time difference quotient, and then apply the quantitative stability estimate \eqref{eq:pairwise} to compare the difference quotients of $u^\varepsilon$ and $u^\delta$ for two parameters $\varepsilon \leq \delta$.

In the remainder of this section, we fix $\kappa=\kappa(n, p)>1$ large enough so that
\begin{equation*}
    \frac{\kappa\alpha}{2}-\frac{2+\alpha}{\alpha}\geq2 \quad \text{and}\quad
    \kappa>\frac2\alpha,
\end{equation*}
where $\alpha \in (0,1)$ is the exponent chosen in Proposition~\ref{prop:regularized}. Then we set
\begin{equation*}
\mu\coloneqq\frac{4(\kappa+2)}{\alpha}.
\end{equation*}

\begin{lemma}\label{lem:noncritical}
Let $1/2\leq s<s+2d\leq2/3$.
There exist constants $r_0\in(0,1)$ and $C>0$, depending only on $n$ and $p$, such that if
\begin{equation*}
0<r\leq\min\{r_0,d^{1/\kappa}\},
\end{equation*}
then, for every $0<\varepsilon\leq r^{\mu}$, we have
\begin{equation}\label{eq:noncritical}
\widetilde M_\varepsilon(s+d,r)
\leq M_{r^{\mu}}(s+2d)+Cr^2.
\end{equation}
\end{lemma}

\begin{proof}
There is nothing to prove if $\widetilde M_\varepsilon(s+d,r)=0$. Otherwise, fix $(x_0,t_0)\in Q_{s+d}$ with $|Du^\varepsilon(x_0,t_0)|\geq r$. Set $\rho=(\theta r)^{1/\alpha}$ for $\theta=\theta(n,p)\in(0,1)$ to be chosen below. Decreasing $r_0$ if necessary, we may assume that $Q_\rho(x_0,t_0)\subset Q_{3/4}$.
By the uniform $C^{1,\alpha}$ estimate \eqref{eq:gradient}, we can fix $\theta$ sufficiently small so that
\begin{equation*}
|Du^\varepsilon(x,t)-Du^\varepsilon(x_0,t_0)| \leq C_0\rho^\alpha=C_0\theta r \leq r/2 \quad\text{for }(x,t)\in Q_\rho(x_0,t_0).
\end{equation*}
In particular, $|Du^\varepsilon|\geq r/2$ in this cylinder, and hence
\begin{equation*}
\|D_q a_\varepsilon(Du^\varepsilon)\|\leq\frac{C}{\sqrt{\varepsilon^2+|Du^\varepsilon|^2}}\leq\frac Cr.
\end{equation*}
Using \eqref{eq:gradient} again, we obtain
\begin{equation*}
\begin{aligned}
\|a_\varepsilon(Du^\varepsilon(x_1,t_1))-a_\varepsilon(Du^\varepsilon(x_2,t_2))\|
&\leq\frac Cr|Du^\varepsilon(x_1,t_1)-Du^\varepsilon(x_2,t_2)|\\
&\leq\frac Cr\left(|x_1-x_2|^\alpha+|t_1-t_2|^{\alpha/2}\right)
\end{aligned}
\end{equation*}
for all $(x_1,t_1),(x_2,t_2)\in Q_\rho(x_0,t_0)$. In other words, the coefficients $(a^{ij}_{\varepsilon}(Du^{\varepsilon}))$ are uniformly elliptic and satisfy
\begin{equation*}
\rho^{\alpha}[a_\varepsilon(Du^\varepsilon)]_{C^\alpha(Q_\rho(x_0,t_0))}\leq\frac{C}{r}\rho^{\alpha}=C\theta.
\end{equation*}
Thus, the interior Schauder estimate (see \cite[Theorem~4.9]{Lie96} for instance) gives
\begin{equation*}
[u_t^\varepsilon]_{C^\alpha(Q_{\rho/2}(x_0,t_0))}
\leq C\rho^{-2-\alpha}\|u^\varepsilon\|_{L^\infty(Q_\rho(x_0,t_0))}\leq C\rho^{-2-\alpha},
\end{equation*}
where $C$ is independent of $\varepsilon$, $r$ and $(x_0,t_0)$.
Consequently,
\begin{equation}\label{eq:time-holder}
\begin{aligned}
|u_t^\varepsilon(x_0,t_0)-u_t^\varepsilon(x_0,t_0-h)|&\leq C\rho^{-2-\alpha}h^{\alpha/2}\\
&\leq Cr^{-(2+\alpha)/\alpha}h^{\alpha/2} \quad \text{for $h\in[0,\rho^2/16]$}.
\end{aligned}
\end{equation}

We now set $h=r^\kappa \in (0, d]$ and $\delta=r^\mu$. Since $\kappa>2/\alpha$, we may further decrease $r_0$ so that $h\leq\rho^2/16$. Utilizing \eqref{eq:time-holder} and $\kappa\alpha/2-(2+\alpha)/\alpha\geq2$, we obtain
\begin{equation*}
\begin{aligned}
    \left|u_t^\varepsilon(x_0,t_0)-\frac{u^\varepsilon(x_0,t_0)-u^\varepsilon(x_0,t_0-h)}{h}\right| &\leq \frac{1}{h}\int_0^h|u_t^\varepsilon(x_0,t_0)-u_t^\varepsilon(x_0,t_0-\tau)|\,d\tau\\
    &\leq Cr^{-(2+\alpha)/\alpha}h^{\alpha/2} \\
&\leq Cr^2.
\end{aligned}
\end{equation*}
Since $\varepsilon\leq\delta$ and $\mu\alpha/4-\kappa=2$, the quantitative stability estimate \eqref{eq:pairwise} yields
\begin{equation*}
\left|\frac{u^\varepsilon(x_0,t_0)-u^\varepsilon(x_0,t_0-h)}{h}
-\frac{u^\delta(x_0,t_0)-u^\delta(x_0,t_0-h)}{h}\right|
\leq\frac{C\delta^{\alpha/4}}{h}=Cr^2.
\end{equation*}
Finally, since $(s+2d)^2-(s+d)^2\geq d\geq h$, we have $\{x_0\}\times[t_0-h,t_0]\subset Q_{s+2d}$ and
\begin{equation*}
\left|\frac{u^\delta(x_0,t_0)-u^\delta(x_0,t_0-h)}{h}\right|
\leq \frac1h\int_{t_0-h}^{t_0}|u_t^\delta(x_0,t)|\,dt\leq M_\delta(s+2d).
\end{equation*}
Combining all these inequalities and taking the supremum over $(x_0,t_0)$ proves \eqref{eq:noncritical}.
\end{proof}

\begin{corollary}\label{cor:one-step}
Let $1/2\leq s<s+2d\leq2/3$, and let $r_0$ be as in Lemma~\ref{lem:noncritical}. If
\begin{equation*}
0<r\leq\min\{r_0,d^{1/\kappa}\},
\end{equation*}
then, for every $0<\varepsilon\leq r^\mu$, we have
\begin{equation*}
M_\varepsilon(s)
\leq M_{r^\mu}(s+2d)+Cd^{-2}r\log\left(e+\frac r\varepsilon\right)+Cr^2,
\end{equation*}
where $C>0$ depends only on $n$ and $p$.
\end{corollary}

\begin{proof}
It immediately follows from Lemmas~\ref{lem:small-gradient} and \ref{lem:noncritical}.
\end{proof}

We now combine Corollary~\ref{cor:one-step} with an iteration argument to prove Theorem~\ref{thm:main}.

\begin{proof}[Proof of Theorem~\ref{thm:main}]
By homogeneity, we may assume that $\|u\|_{L^\infty(Q_1)}\leq1$. We use the regularized solutions $u^{\varepsilon}$ constructed in Proposition~\ref{prop:regularized}.

\smallskip\noindent\textit{Step 1: A recursive estimate.}
Choose $l_{\ast}=l_{\ast}(n, p)\geq64$ sufficiently large that
\begin{equation*}
1+4\mu\log l\leq l/2 \quad \text{and} \quad l^{-4}\leq r_0 \quad\text{for all }l\geq l_{\ast}.
\end{equation*}
For $\varepsilon>0$ satisfying $l=\log(e/\varepsilon)>l_{\ast}$, we set
\begin{equation*}
    r=l^{-4},\quad d=l^{-1},\quad\text{and} \quad\delta=r^\mu=l^{-4\mu}.
\end{equation*}
Then it is immediate to check that
\begin{itemize}
    \item $0<r=l^{-4} \leq \min\{r_0, d^{1/\kappa}\}$;

    \item $\varepsilon < \delta=r^{\mu}$, since $\log(e/\delta)=1+4\mu\log l \leq l/2$.
\end{itemize}
Thus, we apply Corollary~\ref{cor:one-step} to find that
\begin{equation}\label{eq:recursive}
\begin{aligned}
    M_\varepsilon(s) &\leq M_{r^\mu}(s+2d)+Cd^{-2}r\log\left(e+\frac r\varepsilon\right)+Cr^2\\
&\leq M_\delta(s+2/l)+C/l,
\end{aligned}
\end{equation}
whenever $1/2\leq s<s+2/l\leq2/3$.

\smallskip\noindent\textit{Step 2: Iteration.}
Fix $\varepsilon\in(0,1]$, and set
\begin{equation*}
\varepsilon_0=\varepsilon,\quad l_0=\log(e/\varepsilon),\quad s_0=1/2.
\end{equation*}
As long as $l_j>l_{\ast}$, we define sequences
\begin{equation*}
\begin{aligned}
\varepsilon_{j+1}&\coloneqq l_j^{-4\mu},\\
l_{j+1}&\coloneqq\log(e/\varepsilon_{j+1})
=1+4\mu\log l_j,\\
s_{j+1}&\coloneqq s_j+2/l_j.
\end{aligned}
\end{equation*}
By the choice of $l_{\ast}$, we have $l_{j+1} \leq l_j/2$ provided that $l_j>l_{\ast}$.
In particular, $l_j\leq2^{-j}l_0$ and there is a first index $N\geq0$ such that $l_N\leq l_{\ast}$. If $N\geq1$, then $l_{N-1}>l_{\ast}$ and so
\begin{equation*}
s_j=\frac{1}{2}+2\sum_{i=0}^{j-1}\frac1{l_i} \leq \frac{1}{2}+\frac{2}{l_{N-1}}\sum_{i=0}^{N-1}2^{-i} \leq \frac{1}{2}+\frac{4}{l_{\ast}}<\frac{9}{16} \quad \text{for $0 \leq j \leq N$}.
\end{equation*}
The same estimate also holds when $N=0$, with empty sums taken to be zero.

Therefore, for each $j=0,\ldots,N-1$, we apply \eqref{eq:recursive} with
\begin{equation*}
\varepsilon=\varepsilon_j,\quad l=l_j,\quad s=s_j,\quad r=l_j^{-4},\quad d=l_j^{-1},\quad \delta=r^\mu=\varepsilon_{j+1}
\end{equation*}
to obtain that
\begin{equation*}
M_{\varepsilon_j}(s_j) \leq M_{\varepsilon_{j+1}}(s_{j+1})+C/l_j \quad \text{for $0\leq j \leq N-1$}.
\end{equation*}
Summing these inequalities yields
\begin{equation}\label{eq:telescoping}
M_\varepsilon(1/2)\leq M_{\varepsilon_N}(s_N)
+C\sum_{j=0}^{N-1}\frac1{l_j}\leq M_{\varepsilon_N}(s_N)+\frac{2C}{l_{\ast}}.
\end{equation}

\smallskip\noindent\textit{Step 3: A uniform bound.}
Applying Lemma~\ref{lem:logarithmic} with $s=9/16$ and $d=1/16$, we obtain
\begin{equation}\label{eq:terminal}
M_{\varepsilon_N}(s_N)
\leq M_{\varepsilon_N}(9/16)
\leq C\log\frac e{\varepsilon_N}
=Cl_N\leq Cl_{\ast}.
\end{equation}
Combining \eqref{eq:telescoping} and \eqref{eq:terminal},
and recalling that $l_{\ast}$ depends only on $n$ and $p$,
we conclude that
\begin{equation*}
\sup_{0<\varepsilon\leq1}\|u_t^\varepsilon\|_{L^\infty(Q_{1/2})}\leq C,
\end{equation*}
which is the desired uniform time derivative estimate. The other results follow from rather standard arguments; see \cite{LLYZ25} for instance.
\end{proof}

We finally prove Corollary~\ref{cor:almost-optimal-gradient} with the aid of \cite[Theorem~1.3]{APR17}.

\begin{proof}[Proof of Corollary~\ref{cor:almost-optimal-gradient}]
By homogeneity, we may assume that $\|u\|_{L^\infty(Q_1)}\leq1$. An application of  Theorem~\ref{thm:main} shows that $\|u_t\|_{L^{\infty}(Q_{3/4})} \leq K=K(n, p)$.

Fix $t\in[-1/4,0]$ and write $v=u(\cdot,t)$. Then it is immediate to check that 
\begin{equation*}
-K\leq\Delta_p^N v\leq K \quad\text{in }B_{3/4}
\end{equation*}
in the viscosity sense. For $\delta\in(0,1)$, let $v_\delta$ solve
\begin{equation*}
\left\{\begin{aligned}
    v_\delta-\delta\Delta_p^N v_\delta&=v &&\text{in }B_{3/4},\\
    v_\delta&=v &&\text{on }\partial B_{3/4}.
\end{aligned}\right.
\end{equation*}
Existence and uniqueness follow from the comparison principle and Perron's method; see \cite[Section~4]{CIL92}. Since $v-\delta K$ and $v+\delta K$ are respectively a subsolution and a supersolution of this problem, the comparison principle gives
\begin{equation*}
\|v_\delta-v\|_{L^\infty(B_{3/4})}\leq\delta K.
\end{equation*}
Thus, we have
\begin{equation*}
-\Delta_p^N v_\delta=f_\delta \coloneqq\frac{v-v_\delta}{\delta},
\quad \text{where $f_\delta\in C(\overline{B}_{3/4})$ and $\|f_\delta\|_{L^\infty(B_{3/4})}\leq K$}.
\end{equation*}
Together with $\|v_\delta\|_{L^\infty(B_{3/4})}\leq1+K$, \cite[Theorem~1.3]{APR17} yields, for every $\beta\in(0,\alpha_p)$,
\begin{equation*}
\|v_\delta\|_{C^{1,\beta}(\overline{B}_{5/8})} \leq C(n,p,\beta).
\end{equation*}
By compactness and the uniform convergence $v_\delta\to v$, this estimate passes to the limit. Since the constants are independent of $t$, we obtain
\begin{equation*}
\sup_{t\in[-1/4,0]} \|u(\cdot,t)\|_{C^{1,\beta}(\overline{B}_{5/8})} \leq C.
\end{equation*}
Finally, combining the spatial $C^{1,\beta}$ estimate with the time Lipschitz bound, a standard interpolation argument (see \cite[Lemma~3.2]{LLYZ25} for instance) gives
\begin{equation*}
\|u\|_{C^{1,\beta}(\overline{Q}_{1/2})}\leq C. \qedhere
\end{equation*}
\end{proof}

\section{The case \texorpdfstring{$p=1$}{p=1}}\label{sec:p-one}
In this section, we prove Theorem~\ref{thm:p-one}. We use the same viscosity convention as in Section~\ref{sec:preliminaries}. For \eqref{eq:mean-curvature}, the upper and lower semicontinuous extensions at a vanishing gradient are
\begin{equation*}
 \operatorname{tr} X-\lambda_{\min}(X)
 \quad\text{and}\quad
 \operatorname{tr} X-\lambda_{\max}(X), \quad \text{respectively}.
\end{equation*}

The level-set mean curvature flow equation satisfies the geometric homogeneity condition of \cite{CGG91}. Therefore, if $v$ is a viscosity solution and $\Theta:\mathbb{R}\to\mathbb{R}$ is continuous and nondecreasing, then $\Theta\circ v$ is again a viscosity solution. This invariance, including at vanishing gradients, is a key ingredient of the viscosity theory for geometric equations; we refer to \cite{CGG91,Gig06}.

\begin{proof}[Proof of Theorem~\ref{thm:p-one}]
Let $t_0=-1/8$. We claim that the function
\begin{equation*}
 \varphi(x,t)\coloneqq|x|^2+2(n-1)(t-t_0)
\end{equation*}
is a viscosity solution to \eqref{eq:mean-curvature}. Indeed, a direct calculation shows that
\begin{equation*}
    D\varphi=2x, \quad D^2\varphi=2I, \quad \varphi_t=2(n-1).
\end{equation*}
In particular, we have
\begin{equation*}
    \Delta \varphi-\frac{\langle D^2\varphi D\varphi, D\varphi \rangle}{|D\varphi|^2}=2n-2=\varphi_t \quad \text{if $x \neq 0$},
\end{equation*}
i.e., the equation holds in a classical sense. At $x=0$, we have $D^2\varphi=2I$, and both semicontinuous extensions have the value $2(n-1)=\varphi_t$.

We next define
\begin{equation*}
 \Theta(s)\coloneqq
 \begin{cases}
  0& \text{if $s\leq0$},\\
  \displaystyle\frac1{\log(e/s)}& \text{if $0<s<1$},\\
  1& \text{if $s\geq1$}.
 \end{cases}
\end{equation*}
Since $\Theta$ is continuous and nondecreasing, \cite[Theorem~5.6]{CGG91} implies that
$u=\Theta\circ\varphi$ is a viscosity solution in $Q_1$. For all sufficiently small $h>0$,
\begin{equation*}
 u(0,t_0)=0 \quad\text{and} \quad u(0,t_0+h)=\frac1{\log\left(e/(2(n-1)h)\right)}.
\end{equation*}
Therefore, we have
\begin{equation*}
 \lim_{h\to 0^+}\frac{|u(0,t_0+h)-u(0,t_0)|}{h^\beta}=\infty\quad\text{for every }\beta>0.
\end{equation*}
Hence, $u$ is not locally $\beta$-H\"older continuous in time at $(0,t_0)$ for any $\beta\in(0,1)$.
\end{proof}

\end{document}